# Funny Problems!

by Florentin Smarandache, Ph. D.
University of New Mexico
Department of Mathematics
200 College Road
Gallup, NM 87301, USA
E-mail: smarand@unm.edu

**Abstract**.
Thirty original and collected problems, puzzles, and paradoxes in mathematics and physics are explained in this paper, taught by the author to the elementary and high school teachers at the University of New Mexico – Gallup in 1997-8 and afterwards. They have a more educational interest because make the students think different!

**Introduction**.
In this paper, we present a set of recreational mathematical problems. For each "solutin" a funny logic is invented in order to give the problem a sense.

**Proposed Funny Problems**.

1) Prove that 2 = 1.

Solution:

  2 pints = 1 quart!

2) A man weighs the following weights on the following dates. How is this possible?

   6/1/70   150 lbs.
   6/3/70    0 lbs.
   6/5/70    25 lbs.
   6/7/70    0 lbs.
   6/9/70   145 lbs.

Solution:
The man is an astronaut who went to the moon and back.

   Outerspace weightlessness: 0 lbs.
   1/6 of Earth gravity, or gravity of the moon: 25 lbs.

3) If you have a couple of threes and divide them in half, why do you end up with 4 pieces?

Solution:

   33 cut in half horizontally will make four pieces.

4) How 70 > 3 = LOVE?

Solution:

Move the characters of 70 > 3 around.

5) 10 - 1 = 0

Solution:
If you have a stick (1) and an egg (0) and you give away the stick (1) you still have the egg (0) left.

6) All monkeys east bananas.
   I eat bananas.
   Therefore, I am a monkey!

7) Twelve minus one is equal to two.

Solution:
12 - 1 = 2 ( take digit 1 from 12).

8) 7 + 7 = 0.

Solution:
Take the sticks from the 7's and rearrange them to form a rectangular zero   .

9) 3 x 2578 = hell

Solution:
Read your calculator upside down: 7734 (the product of the numbers) becomes hell (approximately).

10) An earthworm is cut down the middle. How many halves are there?

Solution:
One, because the other half can still be one whole earthworm.

11) From two false hypotheses get a true statement.

Examples:

a) Grass is edible.            (False)
   Edible things are green.    (False)
   Therefore, grass is green.  (True)

b) All dogs are poodles.       (False)
   Spot is a dog.              (False)
   Thus, spot is a poodle.     (True)

12) How can you add 3 with 3 and get 8?

Solution:
Turn one of the threes around and put them together to make an 8 (approximately).

13) If 10 trees fall down, and no one is around to hear them falling, how many of the trees fall?

Solution:
Ten.

14) When algebraically 1=0?

Solution:
In a null ring, which is a set with only one element and one binary operation. If we take for "+" and for "*" the same operation, we get a commutative unitary ring.
In this case, the unitary element for "*" (which is normally denoted by "1") and the null element, (which is normally denoted by "0") coincide.

15) When is it possible to have: $1 + 1 = 10$?

Solution:
In base 2.

16) Another logic:
How can we have ten divided by two equal to zero?

Solution:
Ten cookies divided by two kids are eaten and nothing remains!

17) You are lost and walking down a road. You want to get to town and know the road leads to town but don't know which direction. You meet two twin boys. You know one boy always tells the truth and one always lies. The boys know the direction to town. You cannot tell the boys apart and can only
ask one question to one boy to find the direction to town. What question should you ask?

Solution:
Ask either boy what the other boy would say is the direction to town. This would be a lie because if you were asking the dishonest boy he would tell you a lie. If you were asking the honest boy he would tell you the truth about what the dishonest boy would say (which would be a lie) so he would give you the wrong direction. Town would then be in the opposite direction.

18) Why are manhole covers round?
You know, the manholes on the streets, is there a reason why they made them round or could they be square or triangular?

Solution:
Manhole covers are round because a circle cannot fall inside of itself. If they were square, triangular or some other shape they could be dropped into the hole, which would be dangerous to traffic.

19) You have eleven lines. How can you move five lines and still have nine?

Solution:

| | | | | | | | | | | |
               <-move-->

to form

|\| | |\| E

20) You have a cannon and two identical cannon balls. You take the cannon to a large open location that is perfectly flat and you adjust the cannon barrel so that it is perfectly level. You load one of the cannon balls into the cannon and you hold the other cannon ball at the same height as the barrel. You fire the cannon and drop the other cannon ball at the same time. Which cannon ball will hit the ground first?

Solution:
 Both cannon balls should hit the ground at the same time, since gravity acts equally on two objects having the same mass. The cannon barrel was leveled and the cannon ball would begin to fall as it moved forward out of the barrel at the same rate as the cannon ball that was dropped by hand. They would hit at the same time but the cannon ball fired from the cannon would hit the ground far away.

21) I am invisible but can be measured. I affect everyone and everything that is anything. I span the universe and change from place to place. What am I?

Solution:
I am "gravity".

22) The Moon rotates at a rate of one rotation to every 27.3 Earth days and revolves around the Earth at a rate of one revolution to every 27.3 Earth days. This seems a strange coincidence. How does it relate to our
perception of the Moon as viewed from Earth?

Solution:
 People on Earth only see one side of the Moon because the same side is always facing us. If you lived on the far side of the Moon you would never see Earth. Man first saw the far or "dark side of the Moon" in the 1960's.

23) A romantic puzzle:

   GEOMETRY IS THE MEASUREMENT OF THE WORLD, THE GEO, THE SAME GEO WE
   PICTURE OR GRAPH IN GEOGRAPHY. THESE ARE EASY AND SENSIBLE. THE ONE
   I COULD NEVER MAKE HEADS OR TAILS OF, THOUGH IS

       TRI ... GON (O) ... METRY

   METRY IS MEASURE AND TRI IS THREE. BUT WHAT THE HECK'S A GON(O) THAT
   ONE HAS TO HAVE THREE OF IT TO METRY?
   EXPONENTIAL SILLINESS . . .

24) What is a hungry man's multiplication factor?

Solution:

   8 x 8

25) Spell out the number  NINE!

Solution:

   | | | | | | | | | | |
     eleven bars!

26) There are two 24 x 24 corrals. In each corral there are 6 steers. The farmer expects to produce a calf from each steer. How many calves will be produced?

Solution:
Zero! (steers can't produce calves.)

27) How would a mathematician measure the intensity of an earthquake on a meter as in the movie Armageddon?

Solution:
 It is impossible to have an earthquake on a meteor!

28) Fifteen hunters went bear hunting. One killed 2 bears. How many bears have One killed?

Solution:
 Two. ("One" is the name of one of the hunters.)

29) w /2 = u. Find a logic for this equality.

Solution:
Double "u" divided by 2 is "u".

30) A Sorites Paradox in Physics.

The Sorites Paradoxes are associated with Eubulides
of Miletus (fourth century B.C.).  Here it is a such example in physics.

      Statement of a Sorites Paradox in physics:

Our visible world is composed of a totality of invisible particles.

      Proof:

   a) An invisible particle does not form a visible object, nor do two invisible particles, three invisible particles, etc.
However, at some point, the collection of invisible particles becomes large enough to form a visible object, but there is apparently no definite point where this occurs.

   b) A similar paradox is developed in an opposite direction.
It is always possible to remove an atom from an object in such a way that what is left is still a visible object.
However, repeating and repeating this process, at some point, the visible object is decomposed so that the left part becomes invisible, but there is no definite point where this occurs.

    Between <A> and <Non-A> there is no clear distinction, no exact frontier.  Where does <A> really end and <Non-A> does really begin?  We extend Zadeh's fuzzy set term to that of "neutrosophic" concept.

This Sorites paradox does not depend on the observer's subjective vision, because it does not matter how good or weak is his/her vision, for the frontier will not be distinctive but vague, imprecise between visible and invisible.

      References:

[1] Chong Hu, "How do you explain the Smarandache Sorites Paradox?", The MAD Scientist Network: Physics,
http://www.madsci.org/posts/archives/970594003.Ph.q.html.

[2] Amber Iler, Staff, Research Scientist, Veridian – ERIM International, Re: How do you explain the Smarandache Sorites Paradox?, The MAD Scientist Network: Physics,
http://www.madsci.org/posts/archives/970594003.Ph.r.html.

[3] Le, Charles T. Le, Book review of "Definitions, Solved and Unsolved Problems, Conjectures, and Theorems in Number Theory and Geometry" (by Florentin Smarandache, Xiquan Publishing House, Phoenix, 2000), in Zentralblatt fuer Mathematik, Berlin, ZBL: 991.48383, 2000.

[4] Smarandache, Florentin, Sorites Paradox in Physics, mss., 1972, <The Florentin Smarandache papers> Special Collection, Arizona State University, Hayden Library, Tempe, AZ 85287-1006, USA.

[5] Smarandache, Florentin, "Invisible Paradox", in "Neutrosophy. / Neutrosophic Probability, Set, and Logic", 105 p., 22-23, American Research Press, Rehoboth, 1998.

[6] Smarandache, Florentin, "Sorites Paradoxes", in "Definitions, Solved and Unsolved Problems, Conjectures, and Theorems in Number Theory and Geometry", edited by M. L. Perez, 86 p., Xiquan Publishing House, Phoenix,
69-70, 2000.

[7] "The Florentin Smarandache papers" Special Collection, Arizona State University, Hayden Library, Tempe, USA.

[8] "The Florentin Smarandache papers" Special Collection, Archives of American Mathematics, University of Texas at Austin, USA.